\documentclass[12pt, leqno]{article}
\usepackage{amsmath}
\usepackage{amssymb}
\def\q \m#1#2{{\raise1pt\hbox{$#1$}\kern-1pt\big/
               \kern-1pt\raise-1pt\hbox{$#2$}}}

\def\H{{\Bbb H}}

\def\R{{\Bbb R}}
\def\Z{{\Bbb Z}}
\def\Q{{\Bbb Q}}
\def\C{{\Bbb C}}
\def\Aut{{\rm Aut}}
\def\wt{{\rm wt}}
\def\End{{\rm End}}
\def \pf {\noindent {\bf Proof :} \,}

\def\ch{{\rm  ch}}

\def\1{{\bf 1}}
\def\o{{\omega}}
\def\l{{\lambda}}
\def\<{\langle}
\def\>{\rangle}
\def\Res{{\rm Res}}
\def\sym{{\rm Sym}}

\newfam\msbfam

\font\twelmsb=msbm10 at 12pt 
\font\sevenmsb=msbm10 at 7pt \font\fivemsb=msbm10 at 5pt

\textfont\msbfam=\twelmsb
\scriptfont\msbfam=\sevenmsb \scriptscriptfont\msbfam=\fivemsb
\newtheorem{thm}{Theorem}[section]

\newtheorem{lem}[thm]{Lemma}
\newtheorem{lemma}[thm]{Lemma}
\newtheorem{prop}[thm]{Proposition}

\newtheorem{conj}[thm]{Conjecture}
\newtheorem{defn}[thm]{Definition}
\newtheorem{rem}[thm]{Remark}

\newcommand{\cC}{{\cal C}}
\newcommand{\tr}{{\rm Tr}}
\newcommand{\g}{{\frak g}}
\newcommand{\h}{{\frak h}}
\newcommand{\W}{{\cal W}}
\newcommand{\cL}{{\cal L}}

\begin{document}

\renewcommand{\theequation}{\thesection.\arabic{equation}}
\setcounter{equation}{0}

\centerline{\Large {\bf Elliptic genus and vertex operator
algebras}} \vskip 7mm \centerline{\bf Chongying
DONG\footnote{Mathematics Department, University of California,
Santa Cruz, CA 95064, U.S.A. (dong@math.ucsc.edu), Partially
supported by NSF grant DMS-9987656 and a research grant from the
Committee on Research, UC Santa Cruz.},
 Kefeng   LIU\footnote{Department of Mathematics, UCLA,CA 90095-1555,
USA (liu@math.ucla.edu). Partially supported by the Sloan
Fellowship and a NSF grant.}, Xiaonan MA\footnote{Centre de
Math\' ematiques, Ecole Polytechnique, 91128 Palaiseau Cedex,
France (ma@math.polytechnique.fr)} }

\vskip 5mm

{\bf Abstract.} We construct bundles of modules of vertex operator
algebras, and prove the rigidity and vanishing theorem for the
Dirac operator on loop space twisted by such bundles.
 This result generalizes many previous results. \\

{\bf 0 Introduction } \vskip 3mm Let $X$ be a compact smooth spin
manifold. The elliptic genus of Landweber-Stong and Ochanine can
be regarded to be  the index of the formal signature operator on
loop space ${\cal L}X$ (see \cite{W}). It is the index of the
following twisted Dirac operator on $X$
$$D\otimes \otimes_{n\geq 1}\sym_{q^n}(TX\otimes_{\R}\C)\otimes S(TX)\otimes
\otimes_{n>0}\Lambda_{q^n} (TX\otimes_{\R}\C)$$ where $S(TX)$ is
the spinor bundle of $X$ and $D$ is the classical Dirac operator.
Here $q$ is a parameter and for a vector bundle $E$

$$\sym_t(E)= 1+t\,E +t^2 \sym^2(E)+\cdots,\ \Lambda_t(E)=1+t\,E +t^2
\Lambda^2(E)+\cdots$$ are respectively the symmetric and wedge
operation of $E$. This elliptic operator should be considered as
infinite number of twisted Dirac operators by taking the $q$
expansion.

 The following elliptic operator
$$D\otimes \otimes_{n\geq 1}\sym_{q^n}(TX\otimes_{\R}\C)\otimes
\otimes_{n>0} \Lambda_{\pm q^{n-1/2}}(TX\otimes_{\R}\C)$$ were
also studied in \cite{W}. It was conjectured in \cite{W} that all
these elliptic operators are rigid, generalizing the famous
vanishing theorem of Atiyah-Hirzebruch for the $\hat{A}$-genus.
There were several rather interesting proofs of these Witten's
conjectures (see \cite{T}, \cite{BT}, \cite{Liu2}, \cite{LMZ1}).
The one relevant to this paper is the proof given in \cite{Liu2},
\cite{Liu4} where the main idea was to use the modular invariance
of affine Kac-Moody characters.

Note that the fibers of the bundles
$$S(TX)\otimes
\otimes_{n>0}\Lambda_{q^n}
(TX\otimes_{\R}\C)\otimes_{n>0}\Lambda_{\pm
q^{n-1/2}}(TX\otimes_{\R}\C)$$ are level one highest weight
integrable module for affine Kac-Moody Lie algebra $D_l^{(1)}$
where $l$ is half the dimension of $X.$ This explains partially
why the modular invariance of characters of level 1 highest weight
integrable modules for affine algebra $D_l^{(1)}$ enters the proof
of the rigidity \cite{Liu2}.

It is well known that the level one highest weight irreducible
module $L(\Lambda_0)$, where $\Lambda_0$ is the fundamental weight
of $D_l^{(1)}$ corresponding to the index 0, is a vertex operator
algebra and the four level one highest weight integrable modules
for the affine algebra $D_l^{(1)}$ are the irreducible modules for
$L(\Lambda_0).$ In our case, the bundles  $S(TX)\otimes
\otimes_{n>0}\Lambda_{q^n} (TX\otimes_{\R}\C),$
$\otimes_{n>0}\Lambda_{\pm q^{n-1/2}}(TX\otimes_{\R}\C)$ are
$L(\Lambda_0)$-bundles in the sense that each fiber is a module of
$L(\Lambda_0).$ In this paper, we construct very general bundles
such that the corresponding twisted Dirac operators is rigid.
Namely, we twist the Dirac operator by a rather general class of
vertex operator algebra bundles and prove the rigidity property of
the corresponding elliptic operators. The main idea in the proof
of rigidity theorem again is to use the modular invariance of
certain trace functions in the theory of vertex operator algebras.

The study of elliptic operators on loop space twisted by a general
vertex operator algebra bundle in this paper is motivated by
attempting to understand the monstrous moonshine \cite{CN}
geometrically. Borcherds proved \cite{B1} the Conway-Norton's
moonshine conjecture for the McKay-Thompson series associated to
the moonshine vertex operator algebra constructed in \cite{FLM}.
But there is still a lot of interest to understand the genus zero
property for the McKay-Thompson series geometrically. For example,
Hirzebruch proposed to realize the $J$ function as an $\hat
A$-genus \cite{H}. We hope that McKay-Thompson series can be
realized as the Monster equivariant elliptic genus of certain
elliptic operator on some orbifold \cite{DLiuM}. This will lead to
the study of elliptic operators on loop space twisted by vertex
operator algebra bundles whose fibers are twisted modules for the
vertex operator algebra. Another motivation for such construction
is to understand the geometric meaning of elliptic cohomology by
using bundles of vertex operator algebras. Some progress has been
made in this direction. Actually our results indicate that
elliptic cohomology should contain certain vertex operator algebra
bundles.

The setting and the argument in this paper give a uniform
treatment for the Dirac operators on loop spaces when $V$ are the
vertex operator algebras associated to the highest weight
integrable modules for affine Kac-Moody algebras.  More
importantly, our results even in the case that $V$ is a lattice
vertex operator algebra are totally new: many lattice vertex
operator algebras are not ``natural'' modules for loop groups or
affine Kac-Moody algebras. As far as we know, this is also the
first interesting application of the algebraic theory of vertex
operator algebras into geometry and topology. The ideas and
results in this paper can be carried out for the orbifold elliptic
genus in the setting of \cite{DLiuM}. In
this case, one uses twisted sectors or twisted modules for vertex
operator algebras instead of modules. The results on modular
invariance of trace functions in orbifold theory obtained in
\cite{DLM3} is also needed.

This paper is organized as follows. In Section 1 we review the
basic facts about vertex operator algebras and present some
results which will be used in Section 2. In particular, we discuss
the modularity of certain trace functions associated to vertex
operator algebras and their representations. We also give a brief
account of lattice vertex operator algebras for reader's
convenience. In Section 2, we begin with a compact Lie group $G$
which acts on the vertex operator algebra as automorphisms. Using
the principal $G$-bundle we construct the vertex operator algebra
bundles on manifolds. We then consider the Dirac operator on the
loop space twisted by the vertex operator algebra bundles and
prove a rigidity theorem.

{\bf Acknowledgments.}
 Part of this work was done while the third author
was visiting  UCSC. He would like to  thank the Mathematics
Department of  UCSC for its hospitality. This project was
initiated during the authors' visit to the Morningside Center of
Mathematics in Beijing.

\section{ \normalsize Vertex operator algebras and modular invariance}

\setcounter{equation}{0}

In this section we present some results about vertex operator
algebras and their graded traces. We recall some details for those
readers who are not very familiar with the theory of vertex
operator algebras.

\subsection{ \normalsize  Vertex operator algebras and modules}

We give the definitions of vertex operator algebras and their
modules in this section (cf. \cite{B86}, \cite{DLM1}, \cite{FLM},
\cite{Z}).

Let $z,$ $z_0,$ $z_1,$ $z_2$ be commuting formal variables. We
shall use the basic generating function
\begin{eqnarray} \label{0a1}
\delta(z)=\sum_{n\in{\Bbb Z}}z^n,
\end{eqnarray}
which is formally the expansion of the $\delta$-function at $z=1.$
The fundamental
 (and elementary)
properties of the $\delta$-function are in \cite{FLM}, \cite{FHL}
and \cite{DL}.

A {\em vertex  operator algebra} is a ${\Z}$-graded vector space:
\begin{eqnarray} \label{0a2}
V=\bigoplus_{n\in{\Z}}V_n; \ \ \ {\rm for}\ \ \ v\in V_n,\ \
n={\rm wt}\,v;
\end{eqnarray}
such that $\dim V_n<\infty$ for all $n\in\Z$ and $V_n=0$ if $n$ is
sufficiently small; equipped with a linear map
 \begin{align}\label{0a3}
& V \to (\mbox{End}\,V)[[z,z^{-1}]] \\
& v\mapsto Y(v,z)=\sum_{n\in{\Z}}v_nz^{-n-1}\ \ \ \  (v_n\in
\mbox{End}\,V)\nonumber
\end{align}
and with two distinguished vectors ${\bf 1}\in V_0,$ $\omega\in
V_2$ satisfying the following conditions for $u, v \in V$:
\begin{align} \label{0a4}
& u_nv=0\ \ \ \ \ {\rm for}\ \  n\ \ {\rm sufficiently\ large};  \\
& Y({\bf 1},z)=1;  \\
& Y(v,z){\bf 1}\in V[[z]]\ \ \ {\rm and}\ \ \ \lim_{z\to
0}Y(v,z){\bf 1}=v;
\end{align}
and there exists a nonnegative integer $n$ depending on $u,v$ such
that
\begin{align} \label{0a7}
&(z_1-z_2)^n[Y(u,z_1), Y(v,z_2)]=0;\\
&[L(m),L(n)]=(m-n)L(m+n)+\frac{1}{12}(m^3-m)\delta_{m+n,0}(\mbox{rank}\,V)
\end{align}
for $m, n\in {\Z},$ where
\begin{eqnarray} \label{0a9}
L(n)=\omega_{n+1}\ \ \ \mbox{for}\ \ \ n\in{\Z}, \ \ \
\mbox{i.e.},\ \ \ Y(\omega,z)=\sum_{n\in{\Z}}L(n)z^{-n-2}
\end{eqnarray}
and
\begin{align} \label{0a10}
&L(0)v=nv=(\mbox{wt}\,v)v \ \ \ \mbox{for}\ \ \ v\in V_n\
(n\in{\Z}); \\
&\frac{d}{dz}Y(v,z)=[L(-1),Y(u,z)]=Y(L(-1)v,z).
\end{align}
This completes the definition. We denote the vertex operator
algebra just defined by $(V,Y,{\bf 1},\omega)$ (or briefly, by
$V$). The series $Y(v,z)$ are called vertex operators.

An {\em automorphism} $g$ of the vertex operator algebra $V$ is a
linear automorphism of $V$ preserving ${\bf 1}$ and $\omega$  such
that the actions of $g$ and $Y(v,z)$ on $V$ are compatible in the
sense that $gY(v,z)g^{-1}=Y(gv,z)$ for $v\in V.$ Then $gV_n\subset
V_n$ for $n\in\Z.$ The group of all automorphisms of the vertex
operator algebra $V$ is denoted by $\Aut(V)$.

A relevant concept is derivation. A {\em derivation} of $V$ is an
endomorphism $D:V\longrightarrow V$ such that
\begin{eqnarray} \label{0a12}
D(a_nb) = D(a)_nb+a_nD(b)\quad {\rm for\ all\ } a,b\in V, n\in\Z
\end{eqnarray}
and $D(\omega) = 0$. In particular, $D$ preserves the gradation of
$V.$ So the exponential $e^D$ converges on $V$ and is
well-defined. It is easy to see that $e^D$ is an automorphism.

Suppose that $V=\oplus_{n\geq 0}V_n$ with $V_0=\C{\bf 1},$ that
is, $V$ is of CFT type  \cite{DLMM}.  Then  $V_1$ is a Lie algebra
under $[u,v]=u_0v$ with a symmetric invariant bilinear form
$(u,v)=u_1v$ for $u,v\in V_1.$  Moreover,  each $V_n$ is a
$V_1$-module with $u$ acting as $u_0.$ In this case $u_0$ is a
derivation of $V$ and $e^{u_0}$ is an automorphism of $V.$ Set
\begin{eqnarray} \label{0a13}
N=\langle e^{a_0} | a\in V_1\rangle.
\end{eqnarray}
Since $\sigma e^{a_0}\sigma^{-1} = e^{(\sigma a)_0}$ and
$\wt(\sigma(a))= 1$ for any $\sigma\in \Aut(V),$ $N$ is a normal
subgroup of $\Aut(V)$ (cf. \cite{DN}). It is conjectured in
\cite{DN} that $\Aut(V)/N$ is a finite group. An equivalent
conjecture is that the derivation algebra of $V$ is $V_1.$

We say that $V$ is of {\em strong} CFT type if $V$ further
satisfies the condition that $L(1)V_1=0.$ Recall from \cite{FHL}
that a bilinear from $(\cdot, \cdot)$ on $V$ is called {\rm
invariant} if
\begin{eqnarray} \label{0a14}
(Y(u,z)v,w)=(u, Y(e^{zL(1)}(-z^{-2})^{L(0)}v,z^{-1})w)
\end{eqnarray}
for $u,v,w\in V.$ If $V$ is also {\rm simple}, that is $V$ is an
irreducible $V$-module, then there is unique nondegenerate
invariant bilinear from on $V$ \cite{L1}. We shall fix a bilinear
form $(\cdot,\cdot)$ on $V$ so that $(u,v)=u_1v$ for $u,v\in V_1$
(cf. \cite{L1}). It is clear from the definition that
$(gu,gv)=(u,v)$ for any automorphism $g$ and $u,v\in V.$

\begin{rem}\label{r1} {\rm Let $\frak g$ be a finite dimensional simple Lie
algebra and $\hat{\frak g}=\g\otimes\C[t,t^{-1}]\oplus \C K$ the
corresponding affine Kac-Moody Lie algebra. Let $V$ be the vertex
operator algebra associated to the irreducible highest weight
representation of level $m$ for $\hat{\frak g}.$ That is, $V$ is
the irreducible quotient of the Verma module
$U(\hat\g)\otimes_{U(\g\otimes \C[t]+\C K)}\C$ where $\C$ is the
trivial module for $\g\otimes\C[t]$-module and $K$ acts on $\C$ as
scalar $m.$ Then $V_1$ is isomorphic to $\frak g$ and the bilinear
form $(\cdot,\cdot)$ defined by $u_1v$ is the $m$ multiple of
canonical bilinear form $\<\cdot,\cdot\>$ (the square length of a
long root in canonical bilinear form is 2).}
\end{rem}

Now we define admissible modules and ordinary modules for vertex
operator algebras. An admissible $V$-module
$$M=\bigoplus_{n=0}^{\infty}M(n)$$
is a $\Z$-graded vector space with the top level $M(0)\ne 0$
equipped with a linear map
\begin{align}\label{0a15}
&V\longrightarrow (\End\,M)[[z,z^{-1}]]\\
&v\longmapsto\displaystyle{ Y_M(v,z)=\sum_{n\in\Z}v_nz^{-n-1}\ \ \
(v_n\in \End\,M)} \nonumber
\end{align}
which satisfies the following conditions; for $u,v\in V,$ $w\in
M$, $n\in \Z,$
\begin{eqnarray}\label{0a16}\begin{array}{lll}
& &u_nw=0\ \ \
\mbox{for}\ \ \ n\gg 0,\\
& &Y_M({\bf 1},z)=1,\end{array}
\end{eqnarray}
\begin{equation}\label{jacm}
\begin{array}{c}
\displaystyle{z^{-1}_0\delta\left(\frac{z_1-z_2}{z_0}\right)
Y_M(u,z_1)Y_M(v,z_2)-z^{-1}_0\delta\left(\frac{z_2-z_1}{-z_0}\right)
Y_M(v,z_2)Y_M(u,z_1)}\\
\displaystyle{=z_2^{-1}\delta\left(\frac{z_1-z_0}{z_2}\right)
Y_M(Y(u,z_0)v,z_2)}
\end{array}
\end{equation}
(Jacobi identity) where all binomial expressions are to be
expanded in nonnegative integral powers of second variable $z_2$:
This identity is interpreted algebraically as follows: if this
identity is applied to a single vector of $V$ then the coefficient
of each monomial in $z_0,z_1,z_2$ is a finite sum in $V;$
\begin{eqnarray}\label{0a18}
u_mM(n)\subset M(\wt(u)-m-1+n)
\end{eqnarray}
if $u$ is homogeneous. We denote the admissible $V$-module by
$M=(M,Y_M).$

\begin{rem} Let $(M,Y_M)$ be an admissible $V$-module. Then
$L(-1)$-derivation property
\begin{equation}\label{moti}
Y_M(L(-1)v,z)=\frac{d}{dz}Y_M(v,z)
\end{equation}
holds. Moreover, the component operators of $Y_M(\o,z)$ generate a
copy of the Virasoro algebra of central charge ${\rm rank}\,V$
(see \cite{DLM1}).
\end{rem}

A (ordinary)  $V$-module is an admissible $V$-module $M$ which
carries a $\C$-grading induced by the spectrum of $L(0).$ That is,
we have
$$M=\bigoplus_{\lambda \in{\C}}M_{\lambda} $$
where $M_{\l}=\{w\in M|L(0)w=\l w\}.$ Moreover we require that
$\dim M_{\l}$ is finite and for fixed $\l,$ $M_{n+\l}=0$ for all
small enough integers $n.$

A vertex operator algebra $V$ is called {\em rational} if any
admissible $V$-module is a direct sum of irreducible admissible
$V$-modules. It was proved in \cite[Theorem 8.1]{DLM2} (also see
\cite{Z}) that if $V$ is a rational vertex operator algebra then
every irreducible admissible $V$-module is an ordinary $V$-module
and $V$ has only finitely many irreducible admissible  modules up
to isomorphism. We introduce the following notions:\\

(1) $V$ is called {\em holomorphic} if $V$ is rational and $V$ is
the only irreducible module for itself.

(2) $V$ is called $C_2$-{\em cofinite} if $C_2(V)=\<u_2v|u,v\in
V\>$ is of finite codimension.\\

It is a well known conjecture in the theory of vertex operator
algebra that rationality and $C_2$-cofinite conditions are
equivalent. A vertex operator algebra $V$ is said to be {\em
strongly rational} if $V$ is of strong CFT type, rational and
$C_2$-cofinite.

The following theorem was proved in \cite{DM2}.
\begin{thm}\label{red} If $V$ is strongly rational then $V_1$ is a reductive
algebra.
\end{thm}

We shall fix a Cartan subalgebra ${\h}$ of $V_1.$

\subsection{ \normalsize Trace functions and modular invariance}

We first review the vertex operator algebras on torus as defined
in \cite{Z}. The new vertex operator algebra is denoted by $(V,Y[\
],\1, \omega-c/24)$ where $c$ is the central charge of $V.$ The
new vertex operator associated to  a homogeneous element $a$ is
given by
\begin{eqnarray} \label{0a20}
Y[a,z] = \sum_{n\in\Z}a[n]z^{-n-1} = Y(a,e^{z} -1)e^{z\wt{a}}
\end{eqnarray}
while a Virasoro element is $\tilde{\omega} = \omega-c/24$. Thus
\begin{eqnarray} \label{0a21}
a[m] = \Res_z\left(Y(a,z)(\ln{(1+z)})^m(1+z)^{\wt{a}-1}\right)
\end{eqnarray}
and
\begin{eqnarray} \label{0a22}
a[m] = \sum_{i=m}^\infty c(\wt{a},i,m)a(i)
\end{eqnarray}
for some scalars $c(\wt{a},i,m)$ such that $c(\wt{a},m,m)=1.$ In
particular,
\begin{eqnarray} \label{0a23}
a[0]=\sum_{i\geq 0}{{\wt{a}-1}\choose {i}}a(i).
\end{eqnarray}
We also write
\begin{eqnarray} \label{0a24}
L[z] = Y[\omega,z] = \sum_{n\in\Z} L[n]z^{-n-2}.
\end{eqnarray}
Then the $L[n]$ again generate a copy of the Virasoro algebra with
the same central charge $c.$ Now $V$ is graded by the
$L[0]$-eigenvalues, that is
\begin{eqnarray} \label{0a25}
V=\bigoplus_{n\in\Z}V_{[n]}
\end{eqnarray}
where $V_{[n]}=\{v\in V|L[0]v=nv\}.$  It should be pointed out
that for any $n\in\Z$ we have
\begin{eqnarray} \label{0a26}
\sum_{m\leq n}V_n=\sum_{m\leq n}V_{[n]}.
 \end{eqnarray}

It is worthy to remark that if $v\in V_n$ is a lowest weight
vector for the Virasoro algebra generated by $L(m), m\in\Z$ then
$v\in V_{[n]}.$ In particular if $L(1)V_1=0$ then $V_1=V_{[1]}.$

Let $M=\sum_{\lambda\in\C}M_{\lambda}$ be a $V$-module. For
homogeneous $a\in V$ we define
\begin{eqnarray}\label{1a27}
o(a)=a_{\wt a-1},
\end{eqnarray}
and extend $o(a)$ to all $a$ by linearity. Let $a\in V$ we define
\begin{eqnarray}\label{0a27}
Z_M(a,q)=\tr_Mo(a)q^{L(0)-c/24}=q^{-c/24}\sum_{\l\in\C}
(\tr_{M_{\lambda}}o(a))q^{\lambda}.
\end{eqnarray}
If $V$ is $C_2-cofinite$ it is proved in \cite[Theorem 4.4.1]{Z}
that $Z_M(a,q)$ converges to a holomorphic function in upper half
plane with $q=e^{2\pi i\tau}.$

Now we assume that $V$ is rational. Let $M^1,...,M^n$ be the
irreducible $V$-modules. Then there exist rational numbers $\l_i$
for $i=1,...,n$ such that
\begin{eqnarray}\label{new1.28}
M^i=\sum_{p=0}^{\infty}M^i_{\l_i+p}
\end{eqnarray}
(see \cite[Theorem 11.1]{DLM3}) and $M^i_{\l_i}\ne 0.$ For $a\in
V$ we set $Z_i(a,q)=Z_{M^i}(a,q)$ for $i=1,...,n.$  Then
\begin{eqnarray}\label{new1.29}
Z_i(a,q)=q^{\l_i-c/24}\sum_{p=0}^{\infty}(\tr_{M^i_{\l_i+p}}o(a))q^p.
\end{eqnarray}
The following modular property is given in \cite[Theorem
5.5.1]{Z}.
\begin{thm}\label{pa} Assume that $V$ is rational, $C_2$-cofinite.
Let $v\in V_{[m]}$ and $\gamma=\left(\begin{array}{cc}
a & b\\
c & d
\end{array}\right)\in SL(2,\Z).$ Then $Z_s(v,q)$
converges to a holomorphic function in the upper half plane and
there exist scalars ${\gamma}_{st}$ independent of $v,m$ and
$\tau$ such that
\begin{eqnarray}\label{0a28}
Z_s(v,\frac{a\tau+b}{c\tau+d})=(c\tau+d)^m\sum_{t=1}^n{\gamma}_{st}Z_t(v,\tau).
\end{eqnarray}
\end{thm}

We remark that in fact the condition that $V$ is a sum of lowest
weight modules for the Virasoro algebra was assumed in \cite{Z}.
This condition was removed in \cite{DLM3}.

We should also mention that Theorem \ref{pa} does not assert that
$Z_s(v,q)$ is a modular form of weight $m$ on a subgroup of
$SL(2,\Z)$ of finite index. Although the modularity of
$Z_s(v,\tau)$ is always assumed in physics, it is still an open
problem in mathematics. For the discussion of the next section we
introduce the following definition.

\begin{defn}\label{d1} A module $M$ for vertex operator algebra $V$ is
called {\em modular} if there exists a subgroup $\Gamma$ of
$SL(2,\Z)$ of finite index such that
$Z_M(v,\gamma\tau)=\psi(\gamma)(c\tau+d)^mZ_M(v,\tau)$ for $v\in
V_{[m]}$ and $\gamma\in \Gamma$ where $\psi$ is a character on
$\Gamma.$
\end{defn}

The irreducible modules for well-known rational vertex operator
algebras are modular (see \cite{Kac} for the affine vertex
operator algebras and \cite{DMN} for the lattice vertex operator
algebras).

Following \cite{M} we define
\begin{eqnarray}\label{0a29}
Z_s(v,u,q)=\tr_{M^s}e^{2\pi
i(o(v)+(u,v)/2)}q^{L(0)+o(u)+(u,u)/2-c/24}.
\end{eqnarray}
for $u,v\in V_1.$ We remark that the bilinear form on $V_1$ used
in \cite{M} is the minus of the bilinear form used in this paper.
So our $Z_s(v,u,q)$ has a slightly different expression although
they are exactly the same as in \cite{M}. Based on Theorem
\ref{pa}, a modular transformation law is obtained in the
\cite[Main theorem]{M}.

\begin{thm}\label{m} Suppose that $V$ be a rational, $C_2$-cofinite
vertex operator algebra of CFT type. Assume that $u,v\in V_1$ such
that $u,v$ span an abelian Lie subalgebra of $V_1.$ Let
$\gamma=\left(\begin{array}{cc}
a & b\\
c & d
\end{array}\right)\in SL(2,\Z).$ Then $Z_s(u,v,q)$
converges to a holomorphic function in the upper half plane and
\begin{eqnarray}\label{0a30}
Z_s(v,u,\frac{a\tau+b}{c\tau+d})=\sum_{t=1}^n{\gamma}_{st}
Z_t(av+bu,cv+du,\tau)
\end{eqnarray}
where $\gamma_{st}$ is the same as in Theorem \ref{pa}.
\end{thm}
\begin{rem} {\rm Although the modular transformation properties of
 $Z_s(v,u,q)$ was given in \cite{M}, but the convergence of $Z_s(v,u,q)$
was never discussed there. We will prove in the next theorem that
$Z_s(v,0,q)$ is convergent in the upper half plan.}
\end{rem}

{}From now on we assume that $V$ is strongly rational. Recall that
$\h$ is a fixed Cartan subalgebra of the reductive Lie algebra
$V_1$ (cf. Theorem  \ref{red}). Since the homogeneous subspace of
a module $M$  of $V$ is finite dimensional, $M$ is a direct sum of
generalized eigenspaces for $\h.$ Since the restriction of the
bilinear form on $V$ to $\h$ is nondegenerate we can identify
$\h^*$ with $h$ via the bilinear from.

We now define
\begin{eqnarray}\label{chac}
\chi_s(v,\tau)=\chi_{M^s}(v,\tau)= Z_s(v,0,\tau)
\end{eqnarray}
for $v\in V_1.$

\begin{prop}\label{p1.1} Assume that $V$ is strongly rational.
Then $\chi_s(v,\tau)$ converges to a holomorphic function in
${\frak h}\times {\Bbb H}$  where ${\Bbb H}= \{ \tau \in \C; {\rm
Im} \tau >0\}$ is the upper half plane.
\end{prop}

\pf We first recall a result from \cite[Proposition 8]{GN} on the
generators of $V.$ Let $\dim V/C_2(V)=m$ and $x^1,...,x^m$  the
coset representatives such that $L(0)x^i=\mu_ix^i$ and $x^i$ is a
generalized eigenvector for ${\h}$ with eigenvalue $\alpha_i\in
\h.$  Let $U^s$ be the space of lowest weight space of $M^s.$ Then
$M^s=\oplus_{p\geq 0}M^s_{\l_s+p}$ where $M^s_{\l_s}=U^s.$ For
convenience we set $u(m)=u_{\wt u-1-m}$ for $u\in V$ and $m\in
\Z.$ Then $M^s$ is spanned by
$$x^{i_1}(m_1)\cdots x^{i_k}(m_k)U^s$$
for $1\leq i_j\leq m$ and $m_1\geq m_2\geq \cdots \geq m_k>0$ (see
\cite{L3}, \cite{KL}, \cite{GN}, \cite{Bu}). Take $w\in U^s$ to be
a generalized eigenvector for $\h$ with eigenvalue $\gamma\in \h.$
Then the subspace $W$ spanned by
 $$x^{i_1}(m_1)\cdots x^{i_k}(m_k)w$$
for $1\leq i_j\leq m$ and $m_1\geq m_2\geq \cdots \geq m_k>0$ is
invariant under $L(0)$ and $o(v)$ for $v\in\h.$ Since $U^s$ is
finite dimensional, it is enough to prove that $\tr_We^{2\pi
io(v)}q^{L(0)-c/24}$ is holomorphic in $\h\times {\Bbb H}.$

It is easy to see that
$$\tr_We^{2\pi io(v)}q^{L(0)-c/24}\leq q^{-c/24+\l_i}e^{2\pi i(v,\gamma)}
\prod_{j=1}^m\prod_{p>0}(1-q^pe^{2\pi i(v,\alpha_j)})^{-1}$$ where
the inequality holds for each coefficient of the $q^{\mu}e^{2\pi
i}e^{2\pi i(v,t_1\alpha_1+\cdots +t_m\alpha_m)}.$ So it suffices
to show that the power series
$$\prod_{p>0}(1-q^p\xi)^{-1}$$
is convergent absolutely  for $(\xi,\tau)\in \C^* \times \H$.
Write $\xi=e^{\alpha+i\beta}$ and $\tau=x+iy$ for
$\alpha,\beta,x,y\in \R$ with $y>0.$ Then $|q^p\xi|=e^{-2\pi
py}e^{\alpha}.$ It is clear now that $\prod_{p>0}(1-q^p\xi)^{-1}$
is absolutely convergent.
\hfill $\blacksquare$ \\

\begin{rem} {\rm The same argument also shows that $Z_s(u,v,q)$
is convergent for any $u,v\in \h.$}
\end{rem}

Next we discuss the transformation law for $\chi_s(v,\tau)$ under
the modular group action.
\begin{prop}\label{p1} Assume that $V$ is strongly rational.
Then for $\gamma=\left (\begin{array}{l} a \quad b\\
c \quad d
\end{array}  \right ) \in SL(2,\Z),$ $v\in V_1$
\begin{eqnarray}\label{0a30'}
\chi_s({v \over c \tau + d}, { a \tau + b \over c \tau +d}) =
e^{\pi i ( c (v,v)/ (c \tau +
d))}\sum_{j=1}^n\gamma_{sj}\chi_j(v,\tau).
\end{eqnarray}
In particular, if $M$ is modular and $\Gamma$ is the corresponding
subgroup of $SL(2,\Z),$ then
\begin{eqnarray}\label{0a30''}
\chi_M({v \over c \tau + d}, { a \tau + b \over c \tau +d})
=\psi(\gamma)e^{\pi i ( c (v,v)/ (c \tau + d))}\chi_M(v,\tau).
\end{eqnarray}
\end{prop}

\pf It is enough to prove that
\begin{eqnarray}\label{0a31}
\chi_s(\frac{v}{\tau},-\frac{1}{\tau})=e^{\pi
i(v,v)/\tau}\sum_{j=1}^n{S}_{s,t}\chi_t(v,\tau),
\end{eqnarray}
where $S_{s,t}$ corresponds to the matrix
$S=\left(\begin{array}{cc}
0 & -1\\
1 &0
\end{array}\right)$ in Theorem \ref{pa} and
\begin{eqnarray}\label{0a32}
\chi_s(v,\tau+1)=e^{2\pi i(\lambda_s-c/24)}\chi_s(v,\tau).
\end{eqnarray}

The transformation law for $T=\left(\begin{array}{cc}
1 & 1\\
0 & 1
\end{array}\right)$ is clear. For the $S$ matrix we note from Theorem
\ref{m} that
\[
\chi_s(v,-\frac{1}{\tau})=\sum_{t=1}^n{S}_{s,t}Z_t(0,v,\tau)
\]
for any $u\in \h.$ In particular,
\[
\chi_s
(\frac{v}{\tau},-\frac{1}{\tau})=\sum_{t=1}^n{S}_{s,t}Z_t(0,\frac{v}{\tau},\tau).
\]
It is straightforward to verify that
\[
Z_t(0,\frac{v}{\tau},\tau)=e^{\pi i(v,v)/\tau}\chi_t(v,\tau).
\]
This completes the proof.
\hfill $\blacksquare$ \\

\begin{rem}{\rm Recall Remark \ref{r1}. In this case, the factor
$e^{\pi i(v,v)/\tau}$ in Proposition \ref{p1} becomes to $e^{\pi
im\<v,v\>/\tau}.$}
\end{rem}

Recall from Theorem \ref{red} that $V_1$ is a reductive Lie
algebra and ${\frak h}$ is a Cartan subalgebra. Then $V_1={\frak
g}_{ss}\oplus {\frak g}_{a}$ is a direct sum of the semisimple
ideal ${\frak g}_{ss}$ and the center ${\frak g}_{a}$ and ${\frak
h}={\frak h}_{ss}\oplus {\frak g}_{a}$ where ${\frak h}_{ss}$ is a
Cartan subalgebra of ${\frak g}_{ss}.$ Since ${\frak g}_{ss}$ acts
completely on any $V$-module $M^i,$ ${\frak h}_{ss}$ acts on $M^i$
semisimply. For $\alpha\in{\frak h}$ let $M^{i,\alpha}$ be the
generalized eigenspace for $\h$ with eigenvalue $\alpha.$ Also set
$$Q=\{\alpha\in \h| V^{\alpha}\ne 0\}.$$
Then $Q$ is generated by $\alpha_1,...,\alpha_m$ where $\alpha_i$
are given in the proof of Proposition \ref{p1.1}. Clearly
$Q=Q_{ss}\oplus Q_{a}$ where $Q_{ss}$ is a lattice of ${\frak
h}_{ss}$ containing the root lattice of ${\frak g}_{ss}$ and $Q_a$
is a lattice of ${\frak g}_a$ consisting of eigenvalues for
${\frak g}_a$ on $V.$

\begin{lem} There exist $\mu_i\in{\frak h}$  for $i=1,...,n$ such that
\begin{align}\label{0a33}
&V=\oplus_{\alpha\in Q}V^{\alpha},\\
&M^i=\oplus_{\alpha\in Q}M^{i,\mu_i+\alpha}. \nonumber
\end{align}
Moreover, each $V^{\alpha}$ and $M^{i,\mu_i+\alpha}$ are nonzero,
and $Q$ span $\h.$
\end{lem}

\pf Since $V$ is generated by $x^1,...,x^m$ we clearly have
$$V=\oplus_{\alpha\in Q}V^{\alpha}.$$
Showing that $V^{\alpha}$ is nonzero for any $\alpha\in Q$ is
equivalent to showing that if $V^{\beta}$ and $V^{\gamma}$ are
nonzero then $V^{\beta+\gamma}$ and $V^{-\beta}$ are nonzero.
Observe that $u_qV^{\gamma}\subset V^{\beta+\gamma}$ for $u\in
V^{\beta}$ and $q\in\Z.$ It follows from \cite[Proposition
11.9]{DL} that if $u$ is nonzero there there exists $q\in \Z$ such
that $u_qV^{\gamma}$ is nonzero. That is, $V^{\beta+\gamma}$ is
nonzero.

In order to see that $V^{-\beta}$ is nonzero we notice that
$V^{0}\supset V_0$ is nonzero. Since $V$ is simple $V$ is spanned
by $u_qw$ for $u\in V,$ $0\ne w\in V^{\beta}$ and $q\in \Z$ (see
\cite[Corollary 4.2]{DM1} or \cite[Proposition 4.1]{L2}). Thus
$V^{-\beta}$ must be nonzero, otherwise $V^{0}$ would be zero.
Since $M^i$ is an irreducible $V$-module the same argument can
prove that $M^i=\oplus_{\alpha\in Q}M^{i,\mu_i+\alpha}$ and
$M^{i,\mu_i+\alpha}$ is nonzero for any $\alpha\in Q.$

It remains to prove that $Q$ spans $\h.$ Define an invariant
symmetric bilinear form $(\cdot,\cdot)_p$ on $V_1$ such that
$(u,v)_p=\tr_{V_p}o(u)o(v)$ for $u,v\in V.$ It is proved in
\cite{DM2} that if $p$ is big enough the form $(\cdot,\cdot)_p$ is
nondegenerate. Since the form is invariant, the restriction of the
form $(\cdot,\cdot)_p$ to $\h$ is also nondegenerate if $p$ is
large. If $Q$ does not span $\h$ then for some nonzero $u\in \h,$
$o(u)$ has only zero eigenvalue. As a result $(u,v)_p=0$ for all
$v\in \h.$ This is a contradiction.
\hfill $\blacksquare$ \\

This following lemma can be found in \cite{DM2}.

\begin{lem} If $V$ is a strongly rational vertex operator algebra
 then each $V$-module is a completely reducible
$V_1$-module. That is, the action of ${\frak g}_a$ on any
$V$-module is semisimple.
\end{lem}

Let $L$ be the lattice generated by $Q$ and $\mu_i$ for
$i=1,...,n$ and $L^{\circ}=\{\alpha\in\h|(\alpha,L)\subset\Z\}$
the dual lattice of $L.$  Then we immediately have the following
consequence.

\begin{lem} Let $V$ be as before. Then for any $v\in \h,$
$\alpha\in L^{\circ}$ we have
\begin{eqnarray}\label{0a34}
\chi_s(v+\alpha,\tau)=\chi_s(v,\tau).
\end{eqnarray}
\end{lem}

Here is a conjecture on $\chi_M(v,\tau).$
\begin{conj}\label{conj2} {\rm If $V$ is a rational vertex
operator algebra of CFT type then for any irreducible module $M$
\begin{eqnarray}\label{0a35}
\chi_M(v+\alpha\tau,\tau) =e^{-2\pi i(v,\alpha)-\pi
i(\alpha,\alpha)\tau}\chi_M(v,\tau)
\end{eqnarray}
for $v\in \h$ and $\alpha\in L^{\circ}.$}
\end{conj}

This conjecture holds for vertex operator algebras associated to
the highest weight integrable representations for affine Kac-Moody
Lie algebra and for lattice vertex operator algebras. The complete
reducibility of each $V$-module as a $V_1$-module is automatic by
construction in these cases. The transformation property
(\ref{0a35}) for affine vertex operator algebras can be found in
\cite{Kac}. The transformation property (\ref{0a35}) for lattice
vertex operator algebras is discussed below.

Here we briefly recall the structure of lattice vertex operator
algebra $V_K$ associated to a positive definite even lattice $K$
from \cite{B86}, \cite{FLM}, \cite{D}.  Then $K$ is a free abelian
group of finite rank with a $\Z$-valued positive definite
symmetric bilinear form $(\cdot,\cdot)$ such that the square
length of any element is even.  Set $\frak h=K\otimes_{\Z}\C$ and
extend the bilinear form to $\frak h$ by $\C$-linearity. Regarding
$\h$ as an abelian Lie algebra we consider the affine algebra
$\hat\h=\h\otimes\C[t,t^{-1}]\oplus \C C.$ Let $M(1)$ be the
canonical irreducible $\hat\h$-module such that $C$ acts as 1. For
$h\in \h$ and $n\in \Z$ we set $h(n)=h\otimes t^n.$ Then
$M(1)=\C[h(-n)|h\in\h,n>0]$ as vector space. Let $\C[K]$ be the
group algebra of $K.$ The lattice vertex operator algebra $V_K$ is
defined to be $M(1)\otimes \C[K]$ as a vector space. We refer the
reader to \cite{B86} and \cite{FLM} for the definition of the
vertex operators $Y(u,z).$

We write $e^{\alpha}$ for the basis element of $\C[K]$
corresponding to $\alpha\in K.$ For $v=h_1(-n_1)\cdots
h_k(-n_k)\otimes e^{\alpha}$ ($h_i\in\h,$ $n_i>0,$ $\alpha\in K$)
$L(0)v=(\sum_{i=1}^kn_i+(\alpha,\alpha)/2)v.$ So $(V_K)_1$ is
spanned by
$$\{h(-1),e^{\alpha}|h\in\h, \alpha\in K, (\alpha,\alpha)=2\}.$$
Let $\h(-1)$ be the span of $h(-1)$ for $h\in \h.$ Then $\h(-1)$
is a Cartan subalgebra of $(V_K)_1.$ We identify $\h$ with
$\h(-1)$ in an obvious way.

It is well known that $V_K$ is a strongly rational vertex operator
algebras (cf. \cite[Theorem 2.7]{D}, \cite[Theorem 3.16]{DLM1},
\cite[Proposition 12.5]{DLM3}). Let
$K^{\circ}=\{h\in\h|(h,K)\subset \Z\}$ be the dual lattice of $K.$
Let $K^{\circ}=\cup_{g\in K^{\circ}/K}(K+\beta_g)$ be the coset
decomposition. For $g\in  K^{\circ}/K$ we set
$$V_{K+\beta_g}=M(1)\otimes \C[K+\beta_g].$$
 Then $V_{K+\beta_g}$ for $g\in  K^{\circ}/K$ give a complete
list of irreducible modules for $V_K$ up to isomorphism. Now let
$v=h_1(-n_1)\cdots h_k(-n_k)\otimes e^{\alpha}$ with $\alpha\in
K+\beta_g$ we still have
$L(0)v=(\sum_{i=1}^kn_i+(\alpha,\alpha)/2)v.$ Moreover,
$o(h(-1))=h(0)$ acts on $v$ as $(h,\alpha)$ for $h\in\h.$ Note
that the central charge of $V_K$ is the rank of $K.$ Set
$\chi_g(h,\tau)=\tr_{V_{K+\beta_g}}e^{2\pi ih(0)}q^{L(0)-c/24}.$
We have
$$\chi_g(h,\tau)=\frac{\theta_{K+\beta_g}(h,\tau)}{\eta(q)^c}$$
where
$$\theta_{K+\beta_g}(h,\tau)=\sum_{\alpha\in K+\beta_g}e^{2\pi
i(h,\alpha)}q^{(\alpha,\alpha)/2},$$
$$\eta(q)=q^{1/24}\prod_{p>0}(1-q^p),$$
and $c$  is the rank of $K.$ It is immediate to see that
$$\chi_g(h+\alpha\tau,\tau)=e^{-2\pi i(v,\alpha)-\pi
i(\alpha,\alpha)\tau}\chi_g(v,\tau)$$ for $\alpha\in K^{\circ}.$

\section{ \normalsize Rigidity and vanishing theorems of voa elliptic genera}
\setcounter{equation}{0}

In this section we study certain elliptic operators on loop spaces
twisted by  vertex operator algebra bundles. More precisely, we
begin with an arbitrary strong rational vertex operator algebra
$V$ and a compact Lie group $G$ which acts continuously on $V$ as
automorphisms. We use the principal $G$-bundle $P$ and $V$ to form
the associated sequence of vector bundles $\psi(V,P)$ on the
manifold which is used to define elliptic operators. We then show
that these operators are rigid under certain assumptions on the
vertex operator algebra $V.$ We will call the indices of these
twisted elliptic operators the {\em voa elliptic genera}.

\subsection{ \normalsize Rigidity  theorem of voa elliptic genera}

Let $X$ be a  compact manifold and $\dim X= 2k$. We assume that
the $S^1$ acts on $X$, and $TX$ has an $S^1$-equivariant spin
structure. Let $S(TX) = S^+ (TX) \oplus S^- (TX) $ be the spinor
bundle of $TX$. Let $D^X$ be the Dirac operator on $S(TX)$. If $W$
is an $S^1$-equivariant complex vector bundle on $X$, we will
denote by $D\otimes W$ the twisted Dirac operator of $S(TX)\otimes
W$ (cf. \cite[\S 3.3]{BeGeV}).

Recall that the elliptic operator $D\otimes W$ is called rigid if
the equivariant index
$${\rm Ind}_g (D\otimes W)={\rm Tr}\, g_{| {\rm Ker}\,D\otimes W}
 - {\rm Tr}\, g_{| {\rm Coker}\,D\otimes W}$$
 of $D\otimes W$ is constant with respect to
$g\in S^1$.

For a complex (resp. real)  vector bundle $W$ on $X$, as in the
introduction, we let
\begin{eqnarray}  \label{0b1}\begin{array}{l}
 {\rm Sym}_t (W) = 1 + t W + t^2  {\rm Sym}^2 W + \cdots,\\
\Lambda_t (W) = 1 + tW + t^2 \Lambda^2 W + \cdots,
\end{array}\end{eqnarray}
be the symmetric and respectively the exterior power  operations
of $W$ (resp. $W\otimes_\R \C$) in $K(X)[[t]]$.

Let $V= \oplus_{n\geq 0} V_n$ be a strongly rational vertex
operator algebra. Then the bilinear form $(a,b)=a_1b$ on $V_1$ is
$\Aut(V)$-invariant. Recall the automorphism group $N$ of $V$ from
(\ref{0a13}). Let $G$ be a compact Lie group which is contained in
$N$ and acts continuously on $V$ as automorphisms. Then $G$ acts
on each $V$-module.

Let $P$ be an $S^1$-equivariant principal $G$-bundle on $X$. We
define
\begin{eqnarray}\label{0b2}
\psi (V, P) = \sum_{n\geq 0} (P\times_ {G} V_n) q^n \in K(X)[[q]].
\end{eqnarray}
Here $P\times_ {G} V_n$ is the associated vector bundle
corresponding to the representation of $G$ on $V_n$. More
generally, if $M^\mu=\oplus_{p=0}^{\infty}M^{\mu}_{\mu+p}$ is an
irreducible $V$-module, we define
\begin{eqnarray}\label{0b3}
\psi (M^\mu, P) = \sum_{\lambda } (P\times_ {G} M^\mu_\lambda )
q^\lambda \in K(X)[[q^\Q]].
\end{eqnarray}

Recall that the equivariant cohomology group $H^*_{S^1} (X, \Q)$
of $X$ is defined by
\begin{eqnarray}
H^*_{S^1} (X, \Q)= H^*(X \times_{S^1} ES^1, \Q).\nonumber
\end{eqnarray}
where $ES^1$ is the usual  $S^1$-principal bundle over the
classifying space $BS^1$. So $H^*_{S^1} (X, \Q)$ is a module over
$H^*(BS^1, \Q)$ induced by the projection $\overline{\pi} :
X\times _{S^1} ES^1\to BS^1$. Let $p_1(TX)_{S^1} \in H^*_{S^1} (X,
\Q)$ be the equivariant
 first Pontrjagin classes of $TX$.
Also recall that
\begin{eqnarray}
H^*(BS^1, \Q)= \Q [[u]] \nonumber
\end{eqnarray}
with $u$ a generator of degree $2$. Then the $G$-invariant
bilinear form $(\quad)_{V_1}$ defines an $S^1$-equivariant
characteristic class $Q(V_1)_{S^1} $    of $P$.

In the rest, we suppose that there exists $l\in \Z$ such that
\begin{eqnarray}\label{0b5}
Q(V_1)_{S^1}- p_1(TX)_{S^1} =l \cdot  \overline{\pi}^* u^2 \quad
{\rm  in} \quad H^*_{S^1} (X, \Q).
\end{eqnarray}
As in \cite{Liu4}, we call $l$ the anomaly to rigidity.

\begin{thm} \label{th2.1} Assume that $V$ is strong,
rational vertex operator algebra and $M$ an irreducible $V$-module
satisfying (\ref{0a35}). If the $G$-principal bundle $P$ satisfies
(\ref{0b5}), then the elliptic operator
$$D^X \bigotimes(\bigotimes _{m=1}^\infty {\rm Sym}_{q^m} (TX)
\otimes \psi (M, P)) $$ is rigid for $l\leq 0.$ Moreover, its
equivariant index  is zero on $S^1$ if $l<0,$ especially, its
index is zero.
\end{thm}

\begin{rem} {\rm   By combining the argument in this paper and
\cite[\S 2]{LM}, we can easily generalize Theorem \ref{th2.1} to
family case, and obtain the rigidity and vanishing theorems at the
equivariant Chern character level \cite[Definition 2.1]{LM} for
the corresponding  fiberwise twisted Dirac operator of a
fibration. }
\end{rem}

\subsection{\normalsize Proof of Theorem \ref{th2.1}}

For $\tau \in \H = \{ \tau \in \C; {\rm Im} \tau >0\}$,
 $q= e^{ 2\pi i \tau}$, $v\in \C$,  let
\begin{eqnarray}\label{0b6}\begin{array}{l}
\theta (v, \tau)=c(q)q^{1/8} 2 \sin (\pi v ) \Pi_{n=1}^\infty (1 -
q^{n} e^{2 \pi i v}) \Pi_{n=1}^\infty (1 - q^{n} e ^{-2 \pi i v})
\end{array}\end{eqnarray}
be the classical Jacobi theta functions \cite{Ch}, where $c(q)=
\Pi_{n=1}^\infty (1 - q^{n} )$. Set
\begin{eqnarray}\label{0b61}
\theta'(0, \tau) =\left. \frac{\partial \theta(v, \tau) }{\partial
v}\right|_{v=0}.
\end{eqnarray}
Recall  that we have the following transformation formulas of
theta-functions \cite{Ch}:
\begin{eqnarray}\label{0b7}\begin{array}{l}
\theta (t+1, \tau) = -\theta (t,\tau), \qquad
\theta ( t+ \tau, \tau)= - q^{-1/2} e^{- 2 \pi i t} \theta (t, \tau),\\
\theta ({t \over \tau}, - {1 \over \tau})= {1 \over i} \sqrt{\tau
\over i}
 e^{\pi i t^2 \over \tau} \theta (t,\tau), \quad
 \theta (t, \tau+1) = e^{ \pi i \over 4} \theta (t, \tau).
\end{array}\end{eqnarray}

Let $g= e^{2 \pi i t}\in S^1$ be a topological generator of $S^1$.
Let $X^g=\{X_{\alpha}\} $ be the fixed submanifold of the circle
action. Let $i_{\alpha} : X_\alpha \to X$ be the natural
immersion. Let $i_{\alpha}^* : H^*_{S^1} (X, \Q) \to H^*_{S^1}
(X_\alpha, \Q)$ denote the induced homomorphism in equivariant
cohomology. We have the following $S^1$-equivariant decomposition
of $TX$
\begin{eqnarray}\label{0b11}
TX_{|X_{\alpha}} = N_1 \oplus \cdots \oplus N_h \oplus
TX_{\alpha},
\end{eqnarray}
Here $N_{\gamma}$ is a complex vector bundle such that $g$ acts on
it by $e^{2 \pi i m_{\gamma} t}$. We denote the Chern roots of
$N_{\gamma}$ by $2 \pi i x_{\gamma} ^j$, and the Chern roots of
$TX_{\alpha} \otimes_{\R } \C$ by $\{ \pm 2 \pi i y'_j\}$. Let
$\dim_{\C} N_\gamma = d_\gamma$, and $\dim X_\alpha = 2 k_\alpha$.

Now, recall that $P$ is an $S^1$ equivariant $G$-principal bundle
on $X$. We assume that $G$ acts on the right on $P$, and $S^1$
acts on the left on $P$. Let $\omega$ be an $S^1$-equivariant
connection form on $P$, it defines a $S^1$-equivariant horizontal
sunbundle $HP$ of $TP$. Let $\Omega= d \omega +{1\over 2} [\omega,
\omega]$ be the curvature of $\omega$, it is a two form on $P$
with values in $\g$. Let $S$ be the basis of ${\rm Lie} (S^1) =\R$
such that $\exp(tS) =\exp(it)$ for $t\in \R$. Let $S_X, S_P$ be
the
 vector field on $X, P$ induced by $S$. For example, $(S_X f)(x) = {d \over
d \varepsilon}|_{\varepsilon=0} f(\exp (-\varepsilon S)x)$ for
$x\in X$ and $f$ a  $\cC^\infty$ function on $X$.

Let $W$ be a vector space. Let $\rho : G\to {\rm End} (W)$ be a
representation of $G$, let $\W=P\times _G W$ be the corresponding
associated vector bundle. Then the connection $\o$ induces a
connection $\nabla^{\W}$ on $\W$, and the corresponding curvature
is given by $\rho(\Omega)$ \cite[p25]{BeGeV}. Let $\cL^\W (S)$ be
the Lie action of $S$ on the
 $\cC^\infty$ sections of $\W$ which is defined by
$(\cL^\W (S) s)(x) ={d \over d \varepsilon}|_{\varepsilon=0} \Big
( \exp (\varepsilon S)\Big  (s(\exp (-\varepsilon S)x)\Big ) \Big
)$ for $x\in X$ and $s$ a  $\cC^\infty$ section of $\W$ on $X$.
Then the moment  of $S\in {\rm Lie} (S^1) $ is given by
\begin{eqnarray}
\mu(S) = \cL^\W(S) -\nabla^\W_{S_X}.
\end{eqnarray}
Let $TP/X$ be the relative tangent bundle of the fibration $P\to
X$. Let $P^{HP}$, $P^{TP/X}$  be the projections  from $TP= HP
\oplus {TP/X}$ onto  $HP$, $TP/X$. By \cite[p24]{BeGeV}, we know
for $s\in \cC^\infty (P,W)^G = \cC^\infty (X,\W)$. Here
$\cC^\infty (P,W)^G$ is the $G$-invariant $\cC^\infty$ function on
$P$ with values in $W$,
\begin{multline}
\mu(S)s = \cL^\W (S) s -(P^{HP} S_P)\cdot d s
= (S_P-P^{HP} S_P)\cdot d s \\
= (P^{TP/X} S_P)\cdot s
 = - \rho (\o(S_P))s \in \cC^\infty (P,W)^G.
\end{multline}
Now the equivariant curvature of $W$ corresponding to $S$
\cite[p211]{BeGeV} is $\rho (\Omega -  \o(S_P))$. So the
equivariant Chern character of $\W$  for $g=e^{2 \pi i t}$ is
\begin{eqnarray} \label{2b13}
\ch _{g} (\W) = \tr_{W} e^{\rho ({-1\over 2 \pi i}\Omega -2 \pi  t
\o(S_P))}.
\end{eqnarray}

Now we return to our situation. When restricted to $X_\alpha$ we
can calculate (\ref{2b13}) in the following way: when we consider
$P_{|X_\alpha}$ its restriction on $X_\alpha$, then we can define
$f: P_{|X_\alpha} \times S^1 \to G$ by : for $(p,s)\in
P_{|X_\alpha} \times S^1$, $s\cdot p =p\cdot f_p(s)$. Then for $p$
fixed, $f_p: S^1 \to G$ is a group homeomorphism, and for $h\in
G$, $s\in G$,
\begin{eqnarray}
h^{-1} f_{hp}(s) h = f_p(s).
\end{eqnarray}
Now,  we fix $p_0\in P.$ Then  $f_{p_0}(S^1)$ is contained
in a maximum torus $H$ of $G$. Let $G_1$ be the centralize of
$f_{p_0}(S^1)$ in $G$. We choose the Cartan subalgebra $\h$ of $V_1$
in Section 1.1 such that the Lie
algebra of $H$ is contained in $\h$.

By using parallel transport with respect to $HP$, {}from $p_0
G_1$, we get an $S^1$-equivariant $G_1$-principal bundle $P_1$ on
$X_{\alpha}$ which is a subbundle of $P_{|X_{\alpha}}$ and for
$s\in S^1$,  $f_p(s)\in H$ doesn't depend on $p\in P_1$. In fact,
let $X_t : t\in [0,1] \to X$ be a curve such that $x_0$ is the
projection of $p_0$ to $X$, and let $k_t$ be the lift of $x_t$  on
$P$ along $HP$, then $k_t(p_0 g) = k_t(p_0) g$ for any $g\in G$.
In this way, $P_{|X_\alpha} = P_1\times_{G_1} G$ is induced by the
$G_1$-principal bundle $P_1$ on $X_{\alpha}$, and recall that $H$
is the maximum torus of $G_1$. Let $\o_1$ be the restriction of
$\o$ on $P_1$. Then $i \o_1(S_P)$ is constant on $P_1$ which lies
in the lattice $L^{\circ}$ of $\h$ where the lattice
$L^{\circ}$ is defined in Section 1. In the same way, the
restriction of $\Omega$ on $P_1$ lies in $\g_1$, the Lie algebra
of $G_1$. So for $g=e^{2 \pi i t}\in S^1$,
\begin{eqnarray} \label{2b14}
\ch _g (\W) = \tr_{W} e^{ \rho({-1\over 2 \pi i} \Omega - 2 \pi
\o_1 (S_P) t )}.
\end{eqnarray}

We write $T= i \o_1(S_P)\in \h$. Let $U= {-1\over (2 \pi i)^2}
\Omega$. Now from (\ref{1a27}) and (\ref{2b14}), we know that
\begin{eqnarray}\label{0b13}
i^*_{\alpha} \ch _g(P\times_G M_\lambda) = \tr_{M_\lambda} e^{2
\pi i o(U+Tt)}.
 \end{eqnarray}
Thus, the restriction of the equivariant Chern character of $\psi
(M,P)$ on $X_\alpha$ is
\begin{eqnarray}\label{0b131}
 i^*_{\alpha} \ch _g(\psi (M,P)) = \sum_{\lambda}\tr_{M_\lambda}  e^{2 \pi i
o(U+Tt)} q^{L(0)} = q^{c/24} \chi_M (U+Tt,\tau).
\end{eqnarray}
Compare with (\ref{chac}).

For $g=e^{2 \pi i t}, t\in \R$, and $\tau \in \H$, $q= e^{ 2\pi i
\tau}$, we let
\begin{eqnarray} \label{0b14}
F_{M,P}(t, \tau)=q^{-c/24} {\rm Ind}_g (D^X \bigotimes
_{m=1}^\infty {\rm Sym}_{q^m} (TX-\dim X) \otimes \psi (M,P)) .
\end{eqnarray}

For $f(x)$ a holomorphic function, we denote by $f(y')(TX^g) =
\Pi_j f(y'_j)$, the symmetric polynomial which gives
characteristic class of $TX^g$,
 and similarly for $N_\gamma$.
Using the  Atiyah-Bott-Segal-Singer Lefschetz fixed point formula
\cite{AS2},
 (\ref{0b6}), (\ref{0b13}), we find for $t\in [0,1]\setminus \Q$
\begin{eqnarray}  \label{0b15}
\qquad F_{M,P}(t, \tau)= (2\pi i)^{-k} \sum_\alpha \int_{X_\alpha}
\Big [ \theta' (0, \tau)^k \Big ({2 \pi i y' \over \theta (y',
\tau)}\Big )(TX^g)
 { \chi_M (U+Tt,\tau)
\over \Pi_{\gamma} \theta (x_{\gamma} + m_{\gamma} t, \tau)
(N_{\gamma})}\Big ].
\end{eqnarray}

Considered  as functions of $(t, \tau)$, we can obviously  extend
$F_{M,P}( t, \tau)$ to meromorphic functions on $\C \times \H$.
and holomorphic in $\tau$. the first part of Theorem \ref{th2.1}
is  equivalent to the statement that   $F_{M,P}(t, \tau)$ is
independent of $t$. We will prove
  $F_{M,P}(t, \tau)$  is holomorphic on $\C \times \H$, then
Theorem \ref{th2.1} will be deduced from Lemmas \ref{la2.1}.

\begin{lemma} \label{la2.1}  If $Q(V)_{S^1}- p_1(TX)_{S^1} = l \cdot
\overline{\pi}^* u^2 $, then for $a,b\in 2 \Z$,
\begin{eqnarray}\label{0b16}
F_{M,P}(t+ a\tau +b, \tau)= e^{- \pi i l (a^2 \tau + 2 a t)
}F_{M,P}(t, \tau).
\end{eqnarray}
\end{lemma}

\pf  By (\ref{0b7}), for $a,b\in 2 \Z$, $m\in \Z$,  we have
\begin{eqnarray}\label{0b17}
\theta (x + m(t + a \tau + b), \tau) = e^{- \pi i (2m a x + 2 m^2
a t + m^2 a^2 \tau)} \theta (x + mt , \tau).
\end{eqnarray}
By (\ref{0b5}),
\begin{eqnarray}\label{0b18}
 (U+Tt,U+Tt )_{V_1} -\Big ( \sum_j (y'_j)^2
+ \sum_{\gamma,j} (x_{\gamma}^j + m_{\gamma} t)^2\Big ) = l\cdot
t^2
\end{eqnarray}
where $(\cdot,\cdot)_{V_1}$ is the bilinear from on $V_1.$

This means
\begin{eqnarray}\label{0b19}\begin{array}{l}
(T,T)_{V_1}   -\sum_{\gamma} m_{\gamma}^2 d_{\gamma} =l, \quad
 (T, U)_{V_1} = \sum_{\gamma,j} m_{\gamma} x_{\gamma}^j,\\
(U, U)_{V_1} = \sum_j (y'_j)^2 + \sum_{\gamma,j} (x_{\gamma}^j)^2.
\end{array}\end{eqnarray}

By using (\ref{0a34}), (\ref{0a35}),  (\ref{0b15}), (\ref{0b17})
and (\ref{0b19}),
 we get (\ref{0b16}).\hfill $\blacksquare$\\

Now we will  examine the modular transformation property of
$F_{M,P}(t, \tau)$ under the group $SL_2(\Z)$.

For $A= \left ( \begin{array}{l} a \quad b\\
c \quad d
\end{array} \right ) \in SL_2 (\Z)$,  we define its modular
transformation on $\C \times \H$ by
\begin{eqnarray}\label{0b8}\begin{array}{l}
\displaystyle{ A(t, \tau) = \left ( { t \over c \tau + d}, {a \tau
+ b \over c \tau + d} \right ). }
\end{array}\end{eqnarray}
Then the two generators
$ S = \left (  \begin{array}{l} 0 \quad -1 \\
1 \quad \quad 0
 \end{array} \right ), \
T = \left (  \begin{array}{l} 1 \quad 1 \\
0 \quad 1
 \end{array} \right )$
of $SL_2(\Z)$ act on $\C \times \H$ in the following way:
\begin{eqnarray}
S(t,\tau) = \Big ({t \over \tau}, - {1 \over \tau}\Big ),\ T(t,
\tau) = (t, \tau +1). \nonumber
\end{eqnarray}

\begin{lemma} \label{la2.2} For any $A= \left ( \begin{array}{l} a \quad b\\
c \quad d
\end{array} \right ) \in SL_2 (\Z)$, we have
\begin{eqnarray}\label{0b21}
F_{M,P} (A(t, \tau)) = e^{\pi i l ct^2/(c\tau+d)}(c\tau +d)^k
F_{AM,P} (t, \tau),
\end{eqnarray}
where $AM= \Sigma_{\mu} a_{\mu} M^{\mu}$ is a finite complex
linear combination of the irreducible $V$-modules, and we denote
by
\begin{eqnarray}\label{0b22}\begin{array}{l}
\displaystyle{ F_{AM,P} (t, \tau) = (2 \pi i )^{-k} \theta' (0,
\tau)^k \sum_\mu \sum_{\alpha} a_{\mu} \int_{X_\alpha} \Big [
\Big ({2 \pi i y' \over \theta (y', \tau)}\Big )(TX^g) } \\
\displaystyle{\hspace*{40mm} { \chi_{M^\mu} (U+T t, \tau)\over
\Pi_{\gamma} \theta (x_{\gamma} + m_{\gamma} t,
\tau)(N_{\gamma})}\Big ].  }
\end{array}\end{eqnarray}
 the complex linear combination of the  corresponding  equivariant
  indices.
\end{lemma}

\pf     Set
\begin{eqnarray}\label{0b23}
F(t,\tau) = {\theta' (0, \tau)\over \theta (t,\tau)}.
\end{eqnarray}

By (\ref{0b7}), we get
\begin{eqnarray}\label{0b24}
F(A(t, \tau))= (c \tau + d) e ^{-c \pi i t^2/(c \tau + d)} F((c
\tau + d) t ,\tau).
\end{eqnarray}
By Proposition \ref{p1}, (\ref{0b13}), it is easy to see that on
$X_{\alpha}$,
\begin{eqnarray}\label{0b25}
\chi_{M}(A(U+T t, \tau))=e^{c  \pi i
 (U+Tt,U+Tt)_{V_1}/(c \tau + d)} \chi_{AM} (U+T t, \tau).
\end{eqnarray}
with
\begin{eqnarray}\label{0b26}
\chi_{AM} (U+T t, \tau)= \sum_{\mu} a_{\mu} \chi_{M^\mu} (U+T t,
\tau).
\end{eqnarray}

By using (\ref{0b19}), (\ref{0b15}), (\ref{0b24}), (\ref{0b25}),
we get
\begin{eqnarray}\label{0b27}\begin{array}{l}
\displaystyle{ F_{M,P} ({ t \over c \tau + d}, { a \tau + b \over
c \tau + d}) = (2 \pi i)^{-k} \sum_\alpha \int_{X_\alpha}  \Big
[\Big (2 \pi i y'
F(y', { a \tau + b \over c \tau + d} )\Big  )(TX^g)  }\\
\hspace*{25mm}\displaystyle{
 \Pi_{\gamma} \Big (F(x_{\gamma} + {m_{\gamma} t \over c \tau + d},
{ a \tau + b \over c \tau + d}) (N_{\gamma})\Big ) \chi_{M} (U +
{Tt \over c \tau + d}, { a \tau + b \over c \tau + d})
\Big ] }\\
\displaystyle{ \hspace*{10mm} =e^{\pi i l ct^2/(c\tau+d)}(c \tau +
d)^k (2 \pi i)^{-k} \sum_\alpha  \int_{X_\alpha} \Big [\Big (2 \pi
i y'
F((c \tau + d)y', \tau) \Big  )(TX^g)  }\\
\hspace*{25mm} \displaystyle{\Pi_{\gamma} \Big (F((c \tau +
d)x_{\gamma}+ m_{\gamma} t, \tau) (N_{\gamma})\Big ) \chi_{AM} ((c
\tau + d)U + T t, \tau)\Big ]  }
\end{array}\end{eqnarray}

By  (\ref{0b27}), to prove (\ref{0b22}),
  we only need prove the following equation,
\begin{eqnarray}\label{0b28}
\begin{array}{l}
\displaystyle{ \int_{X_\alpha}  \Big [\Big (2 \pi i y '
F((c \tau + d)y', \tau) \Big  )(TX^g) }\\
\hspace*{10mm}\displaystyle{
 \Pi_{\gamma}
\Big (F((c \tau + d)x_{\gamma}+ m_{\gamma} t, \tau)
(N_{\gamma})\Big )
\chi_{AM}  ((c \tau + d)U+T t, \tau)\Big ]   }\\
\displaystyle{ = \int_{X_\alpha}  \Big [\Big (2 \pi i y'
F(y', \tau) \Big  )(TX^g) }\\
\hspace*{25mm} \displaystyle{ \Pi_{\gamma} \Big (F(x_{\gamma}+
m_{\gamma} t, \tau) (N_{\gamma})\Big ) \chi_{AM}  (U+T t,
\tau)\Big ]   }.
\end{array}\end{eqnarray}
By looking at the degree $2 k_{\alpha}$ part, that is the
$k_{\alpha}$-th homogeneous terms of the polynomials in $x$'s,
$y'$'s and $u$'s, on both sides, we get (\ref{0b27}).
The proof of Lemma \ref{la2.2} is complete.\hfill $\blacksquare$\\

The following lemma is a generalization of \cite[Lemma 2.3]{Liu4},
\begin{lemma}  \label{la2.3}
For any $A\in SL_2(\Z)$, the function $F_{A M,P}(t, \tau)$ is
holomorphic in $(t, \tau)$ for $(t, \tau) \in \R \times \H$.
\end{lemma}

\pf Let $z= e^{2 \pi i t}$, and $N= {\rm max}_{\alpha, \gamma}
|m_\gamma|$. Denote by $D_N\subset \C ^2$ the domain
\begin{eqnarray}\label{0b29}\begin{array}{l}
|q|^{1/N} < |z| < |q|^{-1/N}, 0< |q| < 1.
\end{array}\end{eqnarray}
By  (\ref{0b3}), (\ref{0b15})  and (\ref{0b22}), we know that in
$D_N$, $F_{A M,P}(t, \tau)$ has a convergent Laurent series
expansion of the form
\begin{eqnarray}
\sum_\mu a_\mu q^{-c/24+\lambda_{\mu}} \sum_{j=0}^\infty b^A_{j
\mu} (z) q^j
\end{eqnarray}
where $\lambda_{\mu}$ is a rational number such that
$M^{\mu}=\oplus_{p=0}^{\infty}M^{\mu}_{\lambda_{\mu}+p}$ (see
\ref{new1.28}) and  $\{b^A_{j \mu} (z)\}$ are rational functions
of $z$ with possible poles on the unit circle.

Now considered as a formal power series of $q$,
$$\bigotimes _{n=1}^{\infty} {\sym}_{q^n} (TX-\dim X) \otimes \sum_\mu
a_\mu q^{-c/24} \psi (M^{\mu}, P) \Big ) = \sum_\mu a_{\mu}
q^{-c/24+\lambda_{\mu} }\bigotimes _{j=0}^{\infty} V^A_{j,
\mu}q^j$$ with  $V^A_{j, \mu}\in K_{S^1}(X)$. Note that the  terms
in the above two sums correspond to each other. Now, we apply the
Atiyah-Bott-Segal-Singer
 Lefschetz fixed
point formula  to each $V^A_{j, \mu}$,
 for $t \in \R \setminus \Q$, we get
\begin{eqnarray}\label{0b30}
b^A_{j \mu} (z) = {\rm Ind}_z (D \otimes V^A_{j, \mu}).
\end{eqnarray}
This implies that for  $t \in \R \setminus \Q, z= e^{ 2 \pi i t}$,
\begin{eqnarray}\label{0b31}
b^A_{j \mu} (z) = \sum_{l= -N(j)}^{N(j)} a_{l,j}^{A, \mu} z^l.
\end{eqnarray}
for $N(j)$ some positive integer depending on $j$ and $a_{l,j}^{A,
\mu}\in \R$. Since both sides are analytic functions of $z$, this
equality holds for any $z\in \C$.

On the other hand, by multiplying $F_{A M,P}(t, \tau)$ by $f(z) =
\Pi_{\alpha, \gamma} (1- z^{m_\gamma})^{l' d_\gamma}$ $(l'= \dim
X)$, we get  holomorphic functions which have a convergent power
series expansion of the form $\sum_\mu a_\mu q^{-c/24+\l_{\mu}}
\sum_{j=0}^\infty c_{j\mu}^A (z) q^j$,
 with $\{c^A_{j\mu}(z)\}$ polynomial functions in $D_N$.
Comparing the above two expansions, one gets
\begin{eqnarray}\label{0b32}
c^A_{j\mu}(z)= f(z) b^A_{j\mu}(z)
\end{eqnarray}
for each $j$. So by the Weierstrass preparation theorem, we get
$F_{AM,P}(t,
\tau)$ is holomorphic in $D_N$.\hfill $\blacksquare$\\

\begin{lemma} \label{la2.4}  $F_{M,P}(t,\tau)$ is holomorphic on
$\C\times \H$.
\end{lemma}

\pf  Recall that by Proposition \ref{p1.1}, $\chi_M(v,t)$ is
holomorphic on $\h \times \H$. {}From their expressions, we know
the possible polar divisors of $F_{M,P}$ in $\C \times \H$ are of
the form $t= {m \over n} (c\tau + d)$ with $m, n, c, d$ integers
and $(c,d)=1$ or $c=1$ and $d=0$.

We can always find integers $a, b$ such that $ad-bc = 1$, and
consider the matrix
$A=\left ( \begin{array}{l}d\quad -b\\
-c \quad a
\end{array} \right ) \in SL_2(\Z)$. By (\ref{0b21}),
\begin{eqnarray}\label{0b35}\begin{array}{l}
\displaystyle{ F_{A M,P}(t,\tau) = (-c \tau + a)^{-k} F_{M,P} \Big
({t \over -c \tau + a},
 {d \tau -b \over -c \tau + a} \Big )}
\end{array}\end{eqnarray}

Now, if $t= {m \over n}  (c\tau + d)$ is a polar divisor of
$F_{M,P} (t,\tau)$, then one polar divisor of $F_{A M,P}(t,\tau)$
is given by
\begin{eqnarray}\label{0b36}
{t \over -c \tau + a}= {m \over n}
 \Big ( c {d \tau -b \over -c \tau + a} + d \Big ),
\end{eqnarray}
which exactly gives $t = m/n$. This contradicts Lemma \ref{la2.3},
and completes the proof of Lemma  \ref{la2.4}. \hfill $\blacksquare$ \\

{\bf Proof of Theorem \ref{th2.1}}:  By Lemma \ref{la2.4},
$F_{M,P}$ is  holomorphic on $\C\times \H$. For fixed $\tau \in
\H$, if $F_{M,P}(\cdot, \tau)$ isn't identically zero, we let
$\delta$ be the contour $z_0 +2s, z_0 +2 +2s\tau$, $z_0 +2
+2(1-s)\tau$, $z_0 +2(1-s)\tau$ ($s\in [0,1]$), such that
$F_{M,P}(\cdot, \tau)$ does not have any zero on  $\delta$. Then
by (\ref{0b16}),
\begin{eqnarray}\label{0b41}
{1 \over 2\pi i} \int_\delta  {1\over F_{M,P}(t,\tau) } {\partial
\over \partial t} F_{M,P}(t,\tau) dt =4l.
\end{eqnarray}
This means that $F_{M,P}(t,\tau)$ has exactly $4l$ zeros inside
$\delta$.
 Therefore, if $l<0$, $F_{M,P}$ must be identically zero.
If $l=0$, $F_{M,P}(t,\tau)$ is a double periodic holomorphic
function, it  must be independent of $t$. Thus we get  Theorem
\ref{th2.1}.
 \hfill $\blacksquare$ \\

The reader may have noticed that we do {\em not} assume in Theorem
\ref{th2.1} the modularity of $\chi_M(v,\tau).$ When $V$ is the
vertex operator algebra associated to the highest weight
integrable representation for the affine Kac-Moody algebra
$D_l^{(1)}$ it was also proved in \cite[Theorem 5]{Liu4} that
$F_M(v,\tau)$ is a holomorphic Jacobi form by using the fact that
$\chi_M(v,\tau)$ is modular. It turns out the same result holds in
our setting under the assumption that $M$ is modular. See
Definition \ref{d1}.

Recall that a (meromorphic) Jacobi form  of index $m$ and weight
$l$ over $L\rtimes \Gamma$, where $L$ is an integral lattice  in
the complex plane $\C$ preserved by the modular subgroup $\Gamma
\subset SL_2(\Z)$,
 is a (meromorphic) function $F(t, \tau)$ on $\C \times \H$ such that
\begin{eqnarray}\label{0b9}\begin{array}{l}
\displaystyle{F({t \over c \tau + d}, { a \tau + b \over c \tau
+d})
= \psi(A)(c \tau + d)^l e^{2 \pi i m ( c t^2 / (c \tau + d))} F(t, \tau),}\\
\displaystyle{F(t + \lambda \tau + \mu, \tau) = e ^{- 2 \pi i m (
\lambda ^2 \tau + 2 \lambda t)} F(t, \tau),}
\end{array}\end{eqnarray}
where $(\lambda, \mu) \in L$, $A=\left (\begin{array}{l} a \quad b\\
c \quad d
\end{array}  \right ) \in \Gamma,$  and $\psi(A)$ is a character of
$\Gamma.$  If $F$ is holomorphic on $\C \times \H$, we say that
$F$ is a holomorphic Jacobi form.

\begin{thm} \label{th2.2} Let $X, V,M$ and $ P$  be as in Theorem
\ref{th2.1}. Assume that $M$ is modular, then $F_{M,P}$
  is a holomorphic Jacobi form of index $l/2$ and weight $k$ over $(2
  \Z)^2 \rtimes \Gamma$, here $\Gamma$ is the subgroup of $SL(2,\Z)$
  such that $\chi_M(v,\tau)$ is modular over $\Gamma$.
\end{thm}

\pf  Recall that $\Gamma$ is the subgroup of $SL(2,\Z)$ which
defines the modular vertex operator algebra $V$. Then for any
$A=\left ( \begin{array}{l} a \quad b\\
c \quad d \end{array} \right )\in \Gamma$, we have
\begin{eqnarray}\label{0b33}
\chi_{M} (A(U+Tt, \tau)) =\psi(A) e^{ c \pi i (U+Tt, U+Tt)_{V_1}
/(c \tau + d)} \chi_{M} (U+Tt, \tau).
\end{eqnarray}
Now, by (\ref{0b24}), (\ref{0b33}), as Lemma \ref{la2.2}, we get
\begin{eqnarray} \label{0b34}
F_{M,P}(A(t, \tau)) =\psi(A) (c \tau + d)^k e^{\pi i l c t^2/(c
\tau + d)}
    F_{M,P}(t, \tau).
\end{eqnarray}

By Lemmas \ref{la2.1}, \ref{la2.4}, (\ref{0b34}), we get  Theorem
\ref{th2.2}.
\hfill $\blacksquare$\\

\end{document}